# Fast Solver for J2-Perturbed Lambert Problem Using Deep Neural Network

# Fast solver for J2-perturbed Lambert problem using deep neural network

Bin Yang[1] and Shuang Li[2*]

*Nanjing University of Aeronautics and Astronautics, Nanjing 211106, China*

Jinglang Feng[3] and Massimiliano Vasile[4]

*University of Strathclyde, Glasgow, Scotland G1 1XJ, United Kingdom*

**This paper presents a novel and fast solver for the J2-perturbed Lambert problem. The solver consists of an intelligent initial guess generator combined with a differential correction procedure. The intelligent initial guess generator is a deep neural network that is trained to correct the initial velocity vector coming from the solution of the unperturbed Lambert problem. The differential correction module takes the initial guess and uses a forward shooting procedure to further update the initial velocity and exactly meet the terminal conditions. Eight sample forms are analyzed and compared to find the optimum form to train the neural network on the J2-perturbed Lambert problem. The accuracy and performance of this novel approach will be demonstrated on a representative test case: the solution of a multi-revolution J2-perturbed Lambert problem in the Jupiter system. We will compare the performance of the proposed approach against a classical standard shooting**

[1] Ph.D. candidate, Advanced Space Technology Laboratory, No. 29 Yudao Str., Nanjing 211106, China.
[2] Professor, Advanced Space Technology Laboratory, Email: lishuang@nuaa.edu.cn, Corresponding Author.
[3] Assistant Professor, Department of Mechanical and Aerospace Engineering, University of Strathclyde, 75 Montrose Street, Glasgow, UK.
[4] Professor, Department of Mechanical and Aerospace Engineering, University of Strathclyde, 75 Montrose Street, Glasgow, UK.

**method and a homotopy-based perturbed Lambert algorithm. It will be shown that, for a comparable level of accuracy, the proposed method is significantly faster than the other two.**

## I. Introduction

The effect of orbital perturbations, such as those coming from a non-spherical, inhomogeneous gravity field, leads a spacecraft to depart from the trajectory prescribed by the solution of the Lambert problem in a simple two-body model [1], [2]. Since the perturbation due to the J2 zonal harmonics has the most significant effect around all planets in the solar system, a body of research exists that addressed the problem of solving the perturbed Lambert problem accounting for the J2 effect [3], [4]. This body of research can be classified into two categories: indirect methods and shooting methods [5]. Indirect methods transform the perturbed Lambert problem into the solution of a system of parametric nonlinear algebraic equations. For instance, Engles and Junkins [1] proposed an indirect method that uses the Kustaanheimo-Stiefel (KS) transformation to derive a system of two nonlinear algebraic equations. Der [6] presented a superior Lambert algorithm by using the modified iterative method of Laguerre that has good computational performance if given a good initial guess. Armellin et al. [7] proposed two algorithms, based on Differential Algebra, for the multi-revolution perturbed Lambert problems (MRPLP). One uses homotopy over the value of the perturbation and the solution of the unperturbed, or Keplerian, Lambert problem as initial guess. The other uses a high-order Taylor polynomial expansion to map the dependency of the terminal position on the initial velocity, and solves a system of three nonlinear equations. A refinement step is then added to obtain a solution with the required accuracy. A common problem of indirect methods is the need for a good initial guess to solve the system of nonlinear algebraic equations. A bad initial guess increases the time to solve the algebraic system or can lead to a failure of the solution procedure, especially when the transfer time is long.

Shooting methods transcribe the perturbed Lambert problem into the search for the initial velocity vector that provides the desired terminal conditions at a given time. Kraige et al. [8] investigated the efficiency of different shooting approaches and found that a straightforward differential correction algorithm combined with the Rectangular Encke's motion predictor is more efficient than the analytical KS approach. Junkins and Schaub [9] transformed the problem into a two-point boundary value problem and applied Newton iteration method to solve it. The main problem with shooting methods is that, with the increase of the transfer time, the terminal conditions become more sensitive to the variations of the initial velocity and the derivatives of the final states with respect to the initial velocity are more affected by the propagation of numerical errors. In order to mitigate this problem, Arora et al. [10] proposed to compute the derivatives of the initial and final velocity vectors with respect to the initial and final position vectors, and the time of flight, with the state transition matrix. Woollands et al. [11] applied the KS transformation and the modified Chebyshev–Picard iteration to obtain the perturbed solution starting from the solution of the Keplerian Lambert problem, which is to solve the initial velocity vector corresponding to the transfer between two given points with a given time of free flight in a two-body gravitational field [12]. For the multi-revolution perturbed Lambert problem with long flight time, Woollands et al. [13] also utilized the modified Chebyshev-Picard iteration and the method of particular solutions based on the local-linearity, to improve the computational efficiency, but its solution relies on the solution of the Keplerian Lambert problem as the initial guesses. Alhulayil et al. [14] proposed a high-order perturbation expansion method that accelerates convergence, compared to conventional first-order Newton's methods, but requires a good initial guess to guarantee convergence. Yang et al. [15] developed a targeting technique using homotopy to reduce the sensitivity of the terminal position errors on the variation of the initial velocity. However, often techniques that improve robustness of convergence by reducing the sensitivity of the terminal conditions on the initial velocity vector, incur in a higher computational cost.

The major problem of both classes of methods can be identified in the need for a judicious initial guess, often better than the simple solution of the Keplerian Lambert problem. To this end, this paper proposes a novel method combining the generation of a first guess with machine learning and a shooting method based on finite-differences. We propose to train a deep neural network (DNN) to generate initial guesses for the solution of the J2-perturbed Lambert problem and which has been a growing interest in the application of machine learning (ML) to space trajectory design [16], [17]. In Ref. [18] one can find a recent survey of the application of ML to spacecraft guidance dynamics and control. Deep neural network is a technology in the field of ML, which has at least one hidden layer and can be trained using a back-propagation algorithm [18]. Sánchez-Sánchez and Izzo [19] used DNNs to achieve online real-time optimal control for precise landing. Li et al. [16] used DNN to estimate the parameters of low-thrust and multi-impulse trajectories in multi-target missions. Zhu and Luo [20] proposed a rapid assessment approach of low-thrust transfer trajectory using a classification multilayer perception and a regression multilayer perception. Song and Gong [21] utilized a DNN to approximate the flight time of the transfer trajectory with solar sail. Cheng et al. [22] adopted the multi-scale deep neural network to achieve real-time on-board trajectory optimization with guaranteed convergence for optimal transfers. However, to the best of our knowledge ML has not yet been applied to improve the solution of the perturbed Lambert problem.

The DNN-based solver proposed in this paper was applied to the design of trajectories in the Jovian system. The strong perturbation induced by the J2 harmonics of the gravity field of Jupiter creates significant differences between the J2-perturbed and Keplerian Lambert solutions, even for a small number of revolutions. Hence Jupiter was chosen to put the proposed DNN-based solver to the test. The performance of the combination of the DNN first guess generation and shooting will be compared against two solvers: one implementing the homotopy method of Yang et al. [15], the other implementing a direct application of Newton method starting from a first guess generated with the solution of the Keplerian Lambert problem. The homotopy method in Ref. [15] was chosen for its simplicity of implementation and robustness also in the case of long transfer times.

The rest of this paper is organized as follows. In Sec. II, the J2-perturbed Lambert problem and the shooting method are presented. Sec. III investigates eight sample forms and their learning features for the DNN. With comparative analysis of the different sample forms and standardization technologies, the optimal sample form for the J2-perturbed Lambert problem is found. The algorithm using the deep neural network and the finite difference-based shooting method is proposed and implemented to solve the J2-perturbed Lambert problem in Sec. IV. Considering Jupiter's J2 perturbation, Sec. V compares the numerical simulation results of the proposed algorithm, the traditional shooting method and the method with homotopy technique. Finally, the conclusions are made in Sec. VI.

## II.     J2-perturbed Lambert Problem

This section presents the dynamical model we used to study the J2-perturbed Lambert problem and the shooting method we implemented to solve it.

### A.  Dynamical modeling with J2 perturbation

The J2 non-spherical term of the gravity field of planets and moons in the solar system induces a significant variation of the orbital parameters of an object orbiting those celestial bodies. Thus, the accurate solution of the Lambert problem [12] needs to account for the J2 perturbation, especially in the case of a multi-revolution transfer. The dynamic equations of an object subject to the effect of J2 can be written, in Cartesian coordinates, in the following form:

$$\begin{cases} \dot{x} = v_x \\ \dot{y} = v_y \\ \dot{z} = v_z \\ \dot{v}_x = -\dfrac{\mu x}{r^3}\left(1+\dfrac{3}{2}J_2\left(\dfrac{R}{r}\right)^2\left(1-5\dfrac{z^2}{r^2}\right)\right) \\ \dot{v}_y = -\dfrac{\mu y}{r^3}\left(1+\dfrac{3}{2}J_2\left(\dfrac{R}{r}\right)^2\left(1-5\dfrac{z^2}{r^2}\right)\right) \\ \dot{v}_z = -\dfrac{\mu z}{r^3}\left(1+\dfrac{3}{2}J_2\left(\dfrac{R}{r}\right)^2\left(3-5\dfrac{z^2}{r^2}\right)\right) \end{cases} \quad (1)$$

where $\mu$, $R$, and $J_2$ represent the gravitational constant, mean equator radius and oblateness of the celestial body, respectively. ($x$, $y$, $z$, $v_x$, $v_y$, $v_z$) is the Cartesian coordinates of the state of the spacecraft, and $r = \sqrt{x^2 + y^2 + z^2}$ is the distance from the spacecraft to the center of the celestial body.

**B. Shooting Method for the J2-perturbed Lambert Problem**

The classical Lambert problem (or Keplerian Lambert problem in the following) considers only an unperturbed two-body dynamics [12]. However, perturbations can induce a significant deviation of the actual trajectory from the solution of the Keplerian Lambert problem. One way to take perturbations into account is to propagate the dynamics in Eqs. (1) and use a standard shooting method for the solution of two-point boundary value problems.

Fig. 1 depicts the problem introduced by orbit perturbations. The solution of the Keplerian Lambert problem, dashed line, provides an initial velocity $v_0$. Because of the dynamics in Eq.(1), the velocity $v_0$ corresponds to a difference $\Delta r_{f0} = r_f - r_{f0}$ between the desired terminal position $r_f$ and the propagated one $r_{f0}$, when the dynamics is integrated forward in time, for a period *tof*, from the initial conditions [$r_0$, $v_0$]. In order to eliminate this error, one can use a shooting method to calculate a velocity $v$ that corrects $v_0$. Fig. 1 shows an example with two subsequent varied velocity vectors $v_i$ and the corresponding terminal conditions.

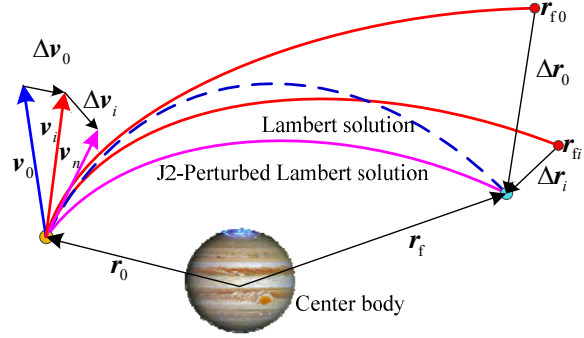

**Fig. 1 Illustration of the shooting method based on Newton's iteration algorithm for the J2-perturbed Lambert problem**

As mentioned in the introduction, shooting methods have been extensively applied to solve the perturbed Lambert problem. Different algorithms have been proposed in the literature to improve both computational efficiency and convergence, e.g. the Picard iteration [11] and the Newton's iteration [23]. In this section, the standard shooting method based on Newton's algorithm is presented [23]. Given the terminal position $\boldsymbol{r}_{fi} = [x_i, y_i, z_i]^T$ and the initial velocity $\boldsymbol{v}_i = [v_{xi}, v_{yi}, v_{zi}]^T$ at the $i$-th iteration, the shooting method requires the Jacobian matrix:

$$\boldsymbol{H}_i = \begin{bmatrix} \dfrac{\partial x_i}{\partial v_{xi}} & \dfrac{\partial x_i}{\partial v_{yi}} & \dfrac{\partial x_i}{\partial v_{zi}} \\ \dfrac{\partial y_i}{\partial v_{xi}} & \dfrac{\partial y_i}{\partial v_{yi}} & \dfrac{\partial y_i}{\partial v_{zi}} \\ \dfrac{\partial z_i}{\partial v_{xi}} & \dfrac{\partial z_i}{\partial v_{yi}} & \dfrac{\partial z_i}{\partial v_{zi}} \end{bmatrix}, \qquad (2)$$

to compute the correction term:

$$\Delta \boldsymbol{v}_i = \boldsymbol{H}^{-1}(\boldsymbol{r}_f - \boldsymbol{r}_i), \qquad (3)$$

where $J^{-1}$ is the inverse of the Jacobian matrix $\boldsymbol{H}_i$, and $\boldsymbol{r}_f$ is the desired terminal position, as shown in Fig. 1. The corrected initial velocity then becomes $\boldsymbol{v}_{i+1} = \boldsymbol{v}_i + \Delta \boldsymbol{v}_i$.

Here the partial derivatives in the Jacobian matrix are approximated with forward differences. Finite differences are computed by introducing a variation $\delta v = 10^{-6}$ in the three components of the initial velocity and computing the corresponding variation of the three components of the terminal conditions $\delta r_{ix}$, $\delta r_{iy}$, and $\delta r_{iz}$. Consequently, the Jacobian matrix can be written as follows.

$$\boldsymbol{H}_i = \left[ \frac{\delta \boldsymbol{r}_{ix}}{\delta v} \quad \frac{\delta \boldsymbol{r}_{iy}}{\delta v} \quad \frac{\delta \boldsymbol{r}_{iz}}{\delta v} \right] \tag{4}$$

Because of the need to compute the Jacobian matrix in Eq. (2), finite-difference-based shooting methods need to perform at least three integrations for each iteration. Furthermore, if the accuracy of the calculation of the Jacobian matrix in Eq.(2) is limited, this algorithm could fail to converge to the specified accuracy or diverge, which is a common situation if the time of flight is long (e.g., tens of revolutions). Homotopy techniques are an effective way to improve the convergence of standard shooting methods for MRPLP but still require an initial guess to initiate the homotopy process and can require the solution of multiple two-point boundary value problems over a number of iterations. Here a DNN is employed to globally map the change in the initial velocity to the variation of the terminal position for a variety of initial state vectors and transfer times. This mapping allows one to generate a first guess for the initial velocity change $\Delta \boldsymbol{v}_i$ by simply passing the required initial state, transfer time and terminal condition as input to the DNN.

In the following, we will present how we trained the DNN to generate good first guesses to initiate a standard shooting method. We will show that an appropriately trained DNN can generate initial guesses that provide improved convergence of the shooting method even for multi-revolution trajectories. It will be shown that the use of this initial guess improves the robustness of convergence of a standard shooting method and makes it significantly faster than the homotopy method in [15].

## III. Sample Learning Feature Analysis

DNN consists of multiple layers of neurons with a specific architecture, which is an analytical mapping from inputs to outputs once its parameters are given. The typical structure of DNN and its neuron computation is illustrated in Fig. 2. The output of each neuron is generated from the input vector $x$, the weights of each component $w$, the offset value $b$, and the activation function $y=f(x)$. The inputs are provided according to the specific problem or the outputs of the neurons of the previous layer. The weight and offset values are obtained through the sample training. The activation function is fixed once the network is built. The training process includes two steps: the forward propagation of the input from the input layer to the output layer; and then the back propagation of the output error from the output layer to the input layer. During this process, the weight and the offset between adjacent layers are adjusted or trained to reduce the error of the outputs.

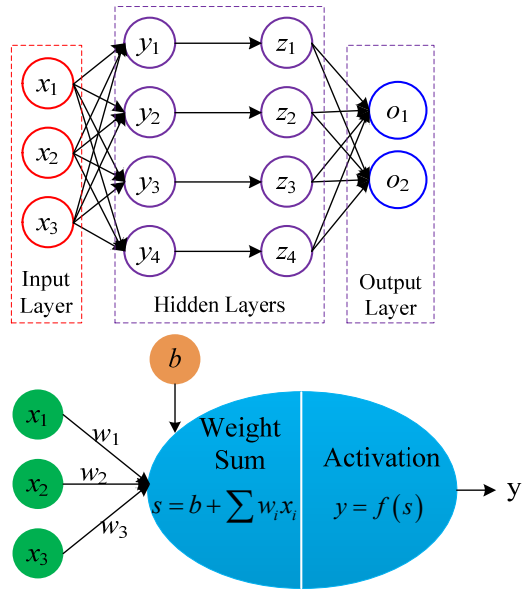

**Fig. 2 The diagram of the DNN structure and neuron computation**

The ability of a DNN to return a good initial guess depends highly on the representation and quality of samples used to train the network. High-quality samples cannot only improve the output accuracy of the network, but also reduce the training cost. Therefore, in the following, we present the procedure used to generate samples with the appropriate features.

A. Definition of Sample Form and Features

In this work two groups of sample forms have been considered: one has the initial velocity $v_0$ solving the J2-perturbed Lambert problem as output, the other has the velocity correction $\Delta v_0$ to an initial guess of $v_0$ as output.

For the first group of sample forms, the input to the neural network includes the known initial and terminal positions $r_0, r_f$ and the time of flight $tof$. The output is only the initial velocity $v_0$ as the terminal velocity can be obtained through orbital propagation once the initial velocity is solved. This type of sample form is defined as

$$S_v = \{[r_0, r_f, tof], [v_0]\} \quad (5)$$

where the subscript 0 and f denotes the start and end of the transfer trajectory, respectively. Thus, when trained with sample form in Eq. (5), the DNN is used to build a functional relationship between $[r_0, r_f, tof]$ and $v_0$.

The second group of sample forms was further divided in two subgroups. One that uses the initial state of the spacecraft $r_0$, the time of flight $tof$ and the terminal error $\Delta r_f$ as input and the other that uses the initial state $r_0$, the time of flight $tof$, the terminal position error $\Delta r_f$ and the initial velocity vector from the Keplerian solution $v_d$ as inputs. These two sample forms are defined as follows:

$$\begin{aligned} S_{dv1} &= \{[r_0, tof, \Delta r_f], [\Delta v_0]\} \\ S_{dv2} &= \{[r_0, v_d, tof, \Delta r_f], [\Delta v_0]\} \end{aligned} \quad (6)$$

In Eq. (6) the output $\Delta v_0$ is always the initial velocity correction $\Delta v_0 = v_0 - v_d$, in which $v_0$ is the initial velocity that solves the J2-perturbed Lambert problem. Thus, when trained with sample forms $S_{dv1}$ and $S_{dv2}$, the DNN realizes a mapping between $\Delta v_0$ and $[r_0, tof, \Delta r_f]$ or $[r_0, v_d, tof, \Delta r_f]$ respectively. The difference between $S_{dv1}$ and $S_{dv2}$ is whether the input includes the initial velocity $v_d$ that is necessary for solving the Jacobian matrix. Therefore, it is theoretically easier to obtain the desired mapping with the input including the initial velocity, i.e. $S_{dv2}$. However, this increases the dimensionality of the sample and might increase the difficulty of training.

For each group of sample forms there are three main ways of parameterizing the state of the spacecraft: Cartesian coordinates, spherical coordinates and the mean orbital elements. Cartesian coordinates provide a general and straightforward way to describe the motion of a spacecraft but state variables change significantly over time even for circular orbits with no orbital perturbations. Spherical coordinates can provide a more contained and simpler variation of the state variables but are singular at the poles. Double averaged mean orbital elements present no variation of semimajor axis, eccentric and inclination due to J2 and a constant variation of argument of the perigee and right ascension of the ascending node [24]. Which parameterization to choose for the training of the DNN will be established in the remainder of this section. The structures of Eqs. (5) and (6) expressed in terms of these three coordinate systems are as follows:

$$\begin{aligned}
S_{v-Car} &= \{[x_0, y_0, z_0, x_f, y_f, z_f, tof], [v_{x0}, v_{y0}, v_{z0}]\} \\
S_{v-Sph} &= \{[r_0, \alpha_0, \beta_0, r_f, \alpha_f, \beta_f, tof], [v_0, \alpha_{v0}, \beta_{v0}]\} \\
S_{v-OEm} &= \{[oe_0, oe_f, tof], [v_0, \alpha_{v0}, \beta_{v0}]\} \\
S_{dv1-Car} &= \{[x_0, y_0, z_0, \Delta x_f, \Delta y_f, \Delta z_f, tof], [\Delta v_{x0}, \Delta v_{y0}, \Delta v_{z0}]\} \\
S_{dv1-Sph} &= \{[r_0, \alpha_0, \beta_0, \Delta r_f, \Delta \alpha_f, \Delta \beta_f, tof], [\Delta v_0, \Delta \alpha_{v0}, \Delta \beta_{v0}]\} \\
S_{dv2-Car} &= \{[x_0, y_0, z_0, v_{xd}, v_{yd}, v_{zd}, \Delta x_f, \Delta y_f, \Delta z_f, tof], [\Delta v_{x0}, \Delta v_{y0}, \Delta v_{z0}]\} \\
S_{dv2-Sph} &= \{[r_0, \alpha_0, \beta_0, v_d, \alpha_d, \beta_d, \Delta r_f, \Delta \alpha_f, \Delta \beta_f, tof], [\Delta v_0, \Delta \alpha_{v0}, \Delta \beta_{v0}]\} \\
S_{dv2-OEm} &= \{[oe_d, \Delta r_f, \Delta \alpha_f, \Delta \beta_f, tof], [\Delta v_0, \Delta \alpha_{v0}, \Delta \beta_{v0}]\}
\end{aligned} \quad (7)$$

where the subscript *Car*, *Sph* and *OEm* denote the Cartesian coordinate, the spherical coordinate and mean orbital elements, respectively. And $x$, $y$, and $z$ are the Cartesian coordinates of the position vector. And $r$, $\alpha$, and $\beta$ are the distance, azimuth, and elevation angle of position vector in the spherical coordinate system. $oe = [a, e, i, \Omega, w, M]^T$ represents the mean orbital elements.

**B. Performance Analysis of Different Sample Forms**

In this section the performance of the eight sample forms defined in Eq.(7) is assessed in order to identify the best one to train the DNN. We always generate a value for the initial conditions starting from an initial set of orbital elements. Values of the orbital parameters for each sample are randomly generated with the *rand* function in MATLAB using a uniform distribution over the intervals defined in Table 1. Note that semimajor axis and eccentricity are derived from the radii of the perijove and apojove. Considering the strong radiation environment of Jupiter and the distribution of Galilean moons, we want to limit the radius of the pericentre $r_p$ of the initial orbit of each sample to be in the interval [$5R_J$, $30R_J$], where $R_J$ = 71492 km is the Jovian mean radius. The value of the inclination is set to range in the interval [0, 1] radians. The time of flight does not exceed one orbital period $T$, which is approximately calculated using the following formula

$$T = 2\pi \sqrt{\frac{a^3}{\mu_J}} \quad (8)$$

where $a$ is the semi-major axis, $a = (r_a + r_p) / 2$.

**Table 1 Parameters' ranges of the sample**

| Parameters | Range |
|---|---|
| Apojove radius $r_a$ (×$R_J$) | [$r_p$, 30] |
| Perijove radius $r_p$ (×$R_J$) | [5, 30] |
| inclination (rad) | [0, 1] |
| RAAN (rad) | [0, $2\pi$) |
| Argument of perigee (rad) | [0, $2\pi$) |
| Mean anomaly (rad) | [0, $2\pi$) |
| *tof* ($T$) | (0, 1) |

The following procedure is proposed to efficiently generate a large number of samples without solving the J2-perturbed Lambert problem:

Step 1: The initial state $[r_0, v_0]$ and time of flight *tof* are randomly generated.

Step 2: The terminal state $[r_f, v_f]$ is obtained by propagating the initial state $[r_0, v_0]$ under the J2 perturbation dynamics model, for the propagation period *tof*.

Step 3: The Keplerian solution $v_d$ is solved from the classical Lambert problem with the initial and terminal position $r_0$, $r_f$ and flight time *tof*.

Step 4: The end state $[r_{fd}, v_{fd}]$ is obtained by propagating the initial Keplerian state $[r_0, v_d]$ under the J2 perturbation dynamics model, and for the propagation period *tof*.

Step 5: The initial velocity correction $\Delta v_0$ and the end state error $\Delta r_f$ are computed with $\Delta v_0 = v_0 - v_d$ and $\Delta r_f = r_f - r_{fd}$.

Using these five steps, we generated 100000 samples and then grouped them in the eight sample forms given in Eq.(7). Before training, a preliminary learning feature analysis is performed on the distribution of sample data and the correlation between the inputs and the output. Specifically, the mean, standard deviation, and magnitude difference coefficients are used to describe the distribution of the data, and the Pearson correlation coefficient is chosen to evaluate the correlation of the data. Their mathematical definitions are given as follows

$$\bar{X} = \frac{\sum_{j=1}^{n} X_j}{n}$$
$$\sigma = \sqrt{\frac{1}{n}\sum_{j=1}^{n}(X_j - \bar{X})^2} \qquad (9)$$
$$\rho = \log\left(\frac{\max(|X|)}{\min(|X| > 0)}\right)$$

where $\bar{X}$ and $\sigma$ are the mean and standard deviation of the data, respectively. And $n$ is the total number of data. $\rho$ denotes the magnitude difference coefficients that assesses the internal diversity of the data.

The statistical characteristics of the variables in the sample are given in Table 2. For the variables described in Cartesian coordinate, the mean values are close to 0 but the standard deviations are generally large. Furthermore, their magnitude difference coefficients are all more than 5, which indicate a large difference in the absolute values of the variables. For the variables described in spherical coordinate, the most of their standard deviations are less than these described in the Cartesian coordinate. In addition, the magnitude difference coefficients of the magnitude of the position and velocity vectors are less than 1. The variables with smaller standard deviation have better performance in the training process. Therefore, the samples with the variables represented in spherical coordinate are easier to learn than those described in Cartesian coordinates.

Table 2 The statistical distributions of the variables in the samples

| Parameters of sample | Mean | Standard deviations | Magnitude difference coefficients |
|---|---|---|---|
| $r_{0\text{-Car}}$ | [-0.014; 0.087; 0.001] | [11.424; 11.424; 1.145] | [5.125; 5.949; 7.077] |
| $r_{0\text{-Sph}}$ | [14.954; 0.007229; 0.000193] | [6.221; 1.815; 0.070] | [0.777; 4.926; 6.607] |
| $r_{f\text{-Car}}$ | [0.001; -0.031; -0.005] | [12.438; 12.469; 1.246] | [5.094; 4.371; 6.442] |
| $r_{f\text{-Sph}}$ | [16.503; -0.005966; -0.000386] | [6.275; 1.813; 0.070] | [0.777; 4.454; 6.321] |
| $v_{0\text{-Car}}$ | [-0.088351; -0.032116; -0.006130] | [8.916; 8.887; 0.895] | [5.450; 5.185; 6.338] |
| $v_{0\text{-Sph}}$ | [12.082577; -0.002034; -0.000529] | [3.647; 1.821; 0.071] | [0.773; 5.257; 5.851] |
| $oe_0$ | [15.771; 0.257895; 0.087045; 3.148016; 3.137227; 3.138286] | [5.225; 0.177; 0.050; 1.813; 1.812; 1.816] | [0.774; 5.273; 4.145; 4.653; 5.422; 4.634] |
| $oe_f$ | [15.771; 0.257850; 0.087045; 3.147935; 3.137600; 3.151625] | [5.225; 0.177; 0.050; 1.813; 1.812; 1.528] | [0.774; 4.721; 4.145; 5.917; 5.228; 5.163] |
| $v_{d\text{-Car}}$ | [-0.087568; -0.031006; -0.006340] | [8.915; 8.886; 0.895] | [5.654; 5.578; 6.673] |
| $v_{d\text{-Sph}}$ | [12.081729; -0.001599; -0.000538] | [3.647; 1.821; 0.071] | [0.774; 5.651; 6.658] |
| $\Delta r_{f\text{-Car}}$ | [-3.162; -11.384; -0.075] | [1369.838; 1395.080; 187.322] | [10.495; 10.828; 11.371] |
| $\Delta \boldsymbol{r}_{\text{f-Sph}}$ | [1154.249; -0.004; -0.001] | [1589.222; 1.817; 0.135] | [10.283; 4.831; 8.576] |
| $oe_{\text{d}}$ | [15.769; 0.258264; 0.087096; 3.147709; 3.136805; 3.139263] | [5.226179; 0.177; 0.051; 1.813; 1.812; 1.814] | [9.556; 4.800; 5.049; 5.240; 5.927; 4.639] |
| $tof$ | 4.023 | 3.220 | 5.481 |
| $\Delta \boldsymbol{v}_{\text{0-Car}}$ | [-0.000782; -0.001109; 0.000210] | [0.326; 0.283; 0.063] | [9.917; 10.432; 10.471] |
| $\Delta \boldsymbol{v}_{\text{0-Sph}}$ | [0.013321; 0.003913; 0.002884] | [0.436; 1.818; 0.503] | [8.948; 4.672; 5.924] |

It is also known that the learning process is easier if the correlation between the input and output of the sample is stronger. Here the Pearson correlation coefficient is used to describe this correlation and is defined as follows

$$R = \frac{\sum_{j=1}^{n}(X_j - \bar{X})(Y_j - \bar{Y})}{\sigma_X \sigma_Y} \tag{10}$$

where $n$ is the total number of sample data. $\bar{Y}$ and $\sigma_Y$ represent the mean and standard deviation of the data $Y$. $\bar{X}$ and $\sigma_X$ denote the mean and standard deviation of the data $X$.

The matrix of the Pearson correlation coefficients of the proposed sample's inputs and outputs are given in Table 3. The elements of Pearson correlation coefficients matrix are the correlation coefficient between the corresponding input and output variables. The signs of the elements indicate positive and negative correlations, respectively. The absolute values of elements represent the strength of correlation. The greater the absolute value is, the stronger the correlation is.

**Table 3 The matrix of the Pearson correlation coefficients of the input and output for different sample forms**

| Sample Forms | Pearson correlation coefficients matrix |
|---|---|
| $S_{v\text{-Car}}$ | $\begin{bmatrix} -0.003 & \mathbf{-0.764} & 0.004 & 0.000 & \mathbf{0.126} & -0.002 & -0.002 \\ \mathbf{0.764} & 0.002 & 0.003 & \mathbf{-0.122} & -0.001 & -0.001 & 0.003 \\ -0.005 & 0.000 & 0.002 & -0.001 & -0.002 & -0.000 & 0.002 \end{bmatrix}$ |
| $S_{v\text{-Sph}}$ | $\begin{bmatrix} \mathbf{-0.898} & 0.005 & 0.001 & \mathbf{-0.459} & -0.003 & 0.000 & \mathbf{-0.344} \\ 0.002 & \mathbf{-0.116} & -0.000 & 0.004 & 0.000 & -0.001 & 0.003 \\ -0.001 & -0.001 & 0.002 & -0.003 & 0.002 & -0.001 & 0.003 \end{bmatrix}$ |
| $S_{v\text{-OEm}}$ | $\begin{bmatrix} \mathbf{-0.587} & \mathbf{0.291} & -0.002 & -0.002 & 0.002 & 0.000 & \mathbf{-0.587} & \mathbf{0.291} & -0.002 & -0.003 & 0.002 & 0.001 & \mathbf{-0.344} \\ 0.002 & 0.000 & 0.000 & 0.002 & 0.003 & 0.007 & 0.002 & 0.000 & 0.000 & 0.002 & 0.003 & 0.001 & 0.003 \\ -0.002 & -0.001 & -0.007 & -0.004 & 0.002 & 0.001 & -0.002 & -0.001 & -0.007 & -0.004 & 0.000 & -0.008 & 0.003 \end{bmatrix}$ |

$S_{dv1\text{-Car}}$ $\begin{bmatrix} -0.006 & \mathbf{-0.011} & 0.002 & \mathbf{-0.049} & \mathbf{-0.053} & \mathbf{-0.046} & 0.002 \\ \mathbf{0.010} & -0.004 & -0.000 & \mathbf{0.037} & \mathbf{-0.041} & \mathbf{-0.013} & -0.002 \\ 0.003 & -0.006 & 0.003 & \mathbf{0.011} & 0.004 & \mathbf{-0.090} & -0.004 \end{bmatrix}$

$S_{dv1\text{-Sph}}$ $\begin{bmatrix} \mathbf{-0.025} & 0.003 & -0.005 & \mathbf{0.081} & 0.004 & -0.002 & \mathbf{0.010} \\ -0.002 & \mathbf{0.377} & 0.000 & 0.002 & \mathbf{0.254} & 0.000 & -0.001 \\ -0.001 & -0.002 & \mathbf{0.512} & 0.001 & 0.001 & \mathbf{0.045} & 0.004 \end{bmatrix}$

$S_{dv2\text{-Car}}$ $\begin{bmatrix} -0.006 & \mathbf{-0.011} & 0.002 & \mathbf{-0.017} & \mathbf{-0.014} & 0.000 & \mathbf{-0.049} & \mathbf{-0.053} & \mathbf{-0.046} & 0.002 \\ \mathbf{0.010} & -0.004 & -0.000 & \mathbf{0.011} & \mathbf{-0.012} & -0.002 & \mathbf{0.037} & \mathbf{-0.041} & \mathbf{-0.013} & -0.002 \\ 0.003 & -0.006 & 0.003 & \mathbf{0.010} & 0.004 & \mathbf{-0.040} & \mathbf{0.011} & 0.004 & \mathbf{-0.090} & -0.004 \end{bmatrix}$

$S_{dv2\text{-Sph}}$ $\begin{bmatrix} \begin{bmatrix} \mathbf{-0.025} & 0.003 & -0.005 \\ -0.002 & \mathbf{0.377} & 0.000 \\ -0.001 & -0.002 & \mathbf{0.512} \end{bmatrix} & \begin{bmatrix} \mathbf{0.032} & 0.004 & -0.005 \\ 0.005 & \mathbf{0.259} & -0.001 \\ 0.003 & 0.004 & \mathbf{0.297} \end{bmatrix} & \begin{bmatrix} \mathbf{0.081} & 0.004 & -0.002 \\ 0.002 & \mathbf{0.254} & 0.000 \\ 0.001 & -0.001 & \mathbf{0.045} \end{bmatrix} & \begin{bmatrix} \mathbf{0.010} \\ -0.001 \\ 0.004 \end{bmatrix} \end{bmatrix}$

$S_{dv2\text{-OEm}}$ $\begin{bmatrix} \mathbf{-0.019} & \mathbf{0.082} & \mathbf{0.076} & -0.001 & -0.004 & 0.001 & \mathbf{0.081} & 0.004 & -0.002 & \mathbf{0.010} \\ 0.002 & -0.003 & 0.002 & 0.004 & 0.002 & 0.001 & 0.002 & \mathbf{0.254} & 0.000 & -0.001 \\ 0.000 & -0.001 & 0.002 & -0.003 & \mathbf{0.010} & -0.002 & 0.001 & -0.001 & \mathbf{0.045} & 0.004 \end{bmatrix}$

First, it is seen that most elements of the matrix are less than 0.01, indicating the correlations between the inputs and the outputs are generally weak. Second, for the first three sample forms of Table 3, the absolute values of all elements for some rows are less than 0.01. This means that some components of the output variable are in weak-correlation with all input variables, and hence the mapping from these output components to the input variables is very difficult to capture. Therefore, samples with the initial velocity as output, i.e. $S_{v\text{-Car}}$, $S_{v\text{-Sph}}$, and $S_{v\text{-OEm}}$, are not deemed to be ideal for the training of the neural network. Third, by comparing the matrix listed in rows 4 to 7 of Table 3, the absolute values of the elements for the samples described in Cartesian coordinates are smaller than those for the samples described in spherical coordinates. Furthermore, for the samples in spherical coordinates, it is seen that the submatrix of each input variable in the Pearson correlation coefficients matrix is a diagonally dominant matrix, where the elements with large absolute values for each input variable are distributed in different rows and columns, and are independent. Therefore, the samples described in the spherical coordinate have better learning features and performance due to the strong correlations. Additionally, for $S_{dv2\text{-Sph}}$ that includes the Keplerian solution $\boldsymbol{v}_d$ as one of the inputs, the correlation with the initial velocity correction $\Delta \boldsymbol{v}_0$ is [**0.032**, 0.004,

-0.005; 0.005, **0.259**, -0.001; 0.003, 0.004, **0.297**], which is diagonally dominant with large diagonal values, which demonstrates that the Keplerian solution is an important input. Finally, for the sample in the mean orbital elements in the last row of Table 3, the matrix only contains a few elements whose absolute values are greater than 0.01, and most of them are distributed in the first row. The mean variations of semimajor axis, eccentricity and inclination are not affected by the J2 perturbation but only by the variation of the initial velocity. Therefore, only the first row in the matrix displays larger values. In addition, the elements in the first six columns of the Pearson correlation matrix of $S_{dv2\text{-}OEm}$ are generally smaller than others in Table 3, because the outputs of the sample is the initial velocity correction, which is calculated using the osculating orbital elements that contain both the long and the short term effects of the J2 perturbation. Thus the correlation using the mean orbital elements is moderate. This would suggest that the sample $S_{dv2\text{-}Sph}$ is the best option for the training of the DNN among the eight tested sample forms. We will now quantify the training performance for each of the eight sample forms by comparing the training convergence of a given DNN. It has to be noted that the structure of the DNN plays a role as well. For example, a high dimensional sample with more variables needs a larger size DNN with more layers and neurons. However, we argue that, since the sample form selection mainly depends on the problem and the dynamics, a better sample form will have better training performance than other sample forms given the same DNN structure. For this reason, it is reasonable to compare sample forms even on DNN structures that are not optimal. The effect of the structure of DNN on the training performance will be discussed in section V.

Some data pretreatment is necessary to facilitate the training process and improve the prediction accuracy. Standardization, normalization and logarithms are used to pre-process data with large ranges or magnitude differences. Tests in this section were performed using a four-layer fully connected DNN with 50 neurons per hidden layer. The activation functions of the hidden layers and the output layer are all Tanh. The Adaptive moment estimation (Adam) [25] was employed for the optimization. The maximum epoch (or number times that the learning algorithm works through the entire training dataset) was set to 10000 and the initial learning rate was set to 0.001. The construction and training of the DNN are based on the Python implementation of TensorFlow. During the training process, the variations of the mean square error (MSE) between the output of the neural network and the output of the sample for different sample forms are given in Fig. 3. The mathematical expression of the MSE is:

$$MSE = \frac{1}{n}\sum_{i=1}^{n}\left(\hat{y}_i - y_i\right)^2 \tag{11}$$

where $n$ is the number of samples, and $\hat{y}_i$ and $y_i$ are the output predicted by the DNN and the true output respectively. Here MSE has no units because data has been normalized before training.

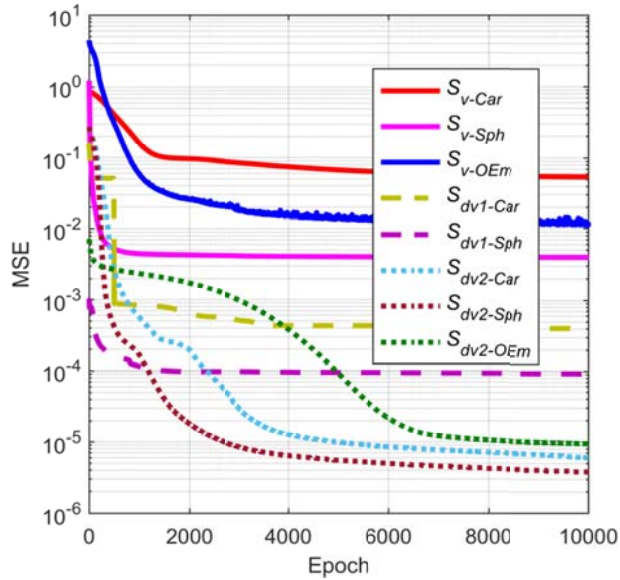

**Fig. 3 The training convergence history for different sample forms**

From Fig. 3 one can see that the MSE of the neural network with the initial velocity correction as the output is significantly smaller than that with the initial velocity as the output. This is because the initial velocity in the sample has a larger range of values and therefore has a more scattered distribution. Also, the MSE of sample $S_{dv1\text{-}Sph}$ is an order of magnitude higher than that of sample $S_{dv2\text{-}Sph}$. Therefore, the accuracy of predicting the initial velocity correction is effectively improved by including the Keplerian velocity in the input of the sample. The blue line in Fig. 3 has obvious fluctuations due to the weak correlations between the output and the input of $S_{v\text{-}OEm}$, as shown in Table 3. Finally, the training results of the samples in spherical coordinate are better than those in Cartesian coordinates, which is consistent with the conclusions drawn in previous sections.

In summary, for the J2-perturbed Lambert problem, the samples described in spherical coordinate appear to be more suitable for the training of a DNN. In fact, among all eight sample forms, the sample form $S_{dv2\text{-}Sph}$ yielded the best learning converge, given the initial position, Keplerian velocity, the terminal position error of the Keplerian solution and time of flight as inputs and the initial velocity correction as output. Therefore, in the remainder of this paper, the $S_{dv2\text{-}Sph}$ sample form is selected for the training of the DNN.

## IV. Solution of the J2-perturbed Lambert Problem Using DNN

The proposed solution algorithm (see the flow diagram in Fig. 4) is made of an Intelligent initial Guess Generator (IGG) and a Shooting Correction Module (SCM). The DNN is used in the IGG to estimate the correction of the Keplerian solution and provide an initial guess to the shooting module. The shooting method discussed in part B of Section II is employed in the SCM to converge to the required accuracy.

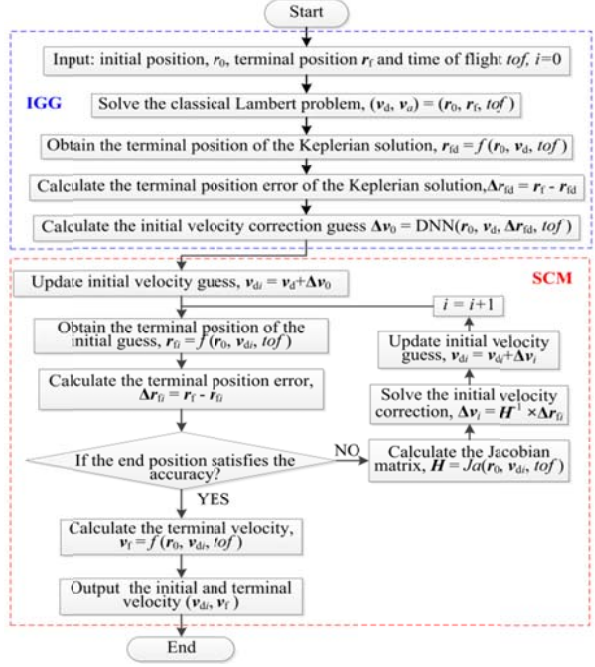

**Fig. 4 The flow chart of the proposed J2-perturbed Lambert problem solver**

As shown in Fig. 4, first the Keplerian Lambert problem is solved with the desired initial and final position vectors. Then the initial conditions $[r_0, v_d]$ are propagated forward in time under the effect of J2 to obtain the terminal position error $\Delta r_{fd}$. With this error, the initial velocity correction is calculated using the trained DNN. The form and the generation method of the samples are described in Section III. Then the finite difference-based shooting method in Section II is applied to correct the initial velocity to make the terminal position meet the rendezvous constraint. The Jacobian matrix is calculated according to Eq.(4), where the partial derivative is approximated with the difference quotient to reduce the computational load.

The method proposed here performs a total of $4i+2$ numerical propagations to obtain the Jacobian matrix and the terminal state, where $i$ is the number of iterations. Additionally, one solution of the Keplerian Lambert problem and one call to the DNN are necessary to obtain the initial velocity guess. Therefore, the calculation time of the proposed method mainly depends on the SCM. As it will be shown in the next section, the initial guess provided by the IGG is close enough to the final solution that the number of iterations required to the SCM to converge to the required accuracy is significantly reduced.

## V.     Case Study of Jupiter Scenario

In this section, taking the Jovian system as an example, some numerical simulations are performed to demonstrate the effectiveness and efficiency of the proposed J2-perturbed Lambert solver. Firstly, different network structures and training parameters are tested to find the optimal ones for this application. Then, we simulate the typical use of the proposed solver with a Monte Carlo simulation whereby a series of transfer trajectories are computed starting from a random set of boundary conditions and transfer times. To be noted that although the tests in this section use the $J_2$, $\mu$ and $R$, of Jupiter the proposed method can be generalized to other celestial bodies by training the corresponding DNNs with a different triplet of values $J_2$, $\mu$ and $R$, but using the same sample form.

### A. DNN Structure Selection and Training

With reference to the results in Section III, the samples used to train the DNN include the initial position, the initial velocity, coming from the solution of the Keplerian Lambert problem, the terminal position error of the Keplerian solution, and the time of flight. The output is the initial velocity correction of the Keplerian solution and all vectors in a sample are expressed in spherical coordinates. In order to generalize the applicability of this method, the ranges of the parameters of the sample given in Table 1 have been appropriately expanded. The range of orbital inclinations is $[0, \pi]$ in radian. The range of times of flight is now in the open interval $(0, 10T)$, where $T$ is calculated using Eq. (8) from the initial state $(\boldsymbol{r}_0, \boldsymbol{v}_0)$. The ranges of other parameters are consistent with Table 1. In total, 200000 training samples are obtained using the rapid sample generation algorithm given in part B of Section III.

Since the structure and training parameters of the neural network also plays a significant impact on the training results, in this section we analyze different DNN structures and settings. Note that once the structure is optimized one would need to loop back and check the optimality of the sample form, however, in this paper we assume that the sample form remains reasonably good even once the DNN structure is changed.

We start by defining the activation functions. Tanh and ReLU are the common activation functions for deep learning while Sigmoid functions are less used because the gradient tends to vanish [26], thus in the following Tanh and ReLU will be used. The output ranges of Tanh and ReLU are [-1, 1] and [0, ∞] respectively, as shown in Fig. 5. The spherical coordinates (magnitude, azimuth, and elevation) of the output of the sample are [0, ∞], [0, $2\pi$] and [-0.5$\pi$, 0.5$\pi$]. Because the range of elevation angle can be transformed from [-0.5$\pi$, 0.5$\pi$] to [0, $\pi$], the ranges of the three components of the spherical coordinates can all meet the requirements of ReLU. Therefore, ReLU is chosen as the activation function of the output layer.

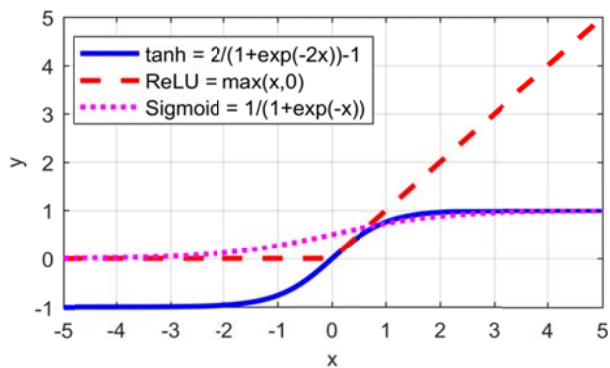

**Fig. 5 The typical activation functions for DNN**

Also in this case the Adaptive moment estimation is used as optimizer. The maximum epoch is 50000 and the other training parameters are the same as in Section III. The training results of DNNs with different sizes are listed in Table 4.

Table 4 Training results of DNNs with different sizes

| Hidden Layers | Neurons per hidden layer | activation function | MSE | Training time (s) |
|---|---|---|---|---|
| 2 | 20 | ReLU | 9.423e-05 | 762 |
| 2 | 20 | Tanh | 3.286e-05 | 839 |
| 2 | 50 | ReLU | 1.435e-05 | 951 |
| 2 | 50 | Tanh | 1.226e-05 | 1084 |
| 2 | 100 | ReLU | 9.423e-06 | 1210 |
| 2 | 100 | Tanh | 9.163e-06 | 1425 |
| 3 | 20 | ReLU | 2.423e-06 | 1198 |
| 3 | 20 | Tanh | 2.154e-06 | 1267 |
| 3 | 50 | ReLU | 1.315e-06 | 1347 |
| 3 | 50 | Tanh | 1.258e-06 | 1523 |
| 3 | 100 | ReLU | 5.631e-06 | 1746 |
| 3 | 100 | Tanh | 1.226e-06 | 1935 |
| 4 | 20 | ReLU | 9.423e-06 | 1648 |
| 4 | 20 | Tanh | 3.286e-06 | 1864 |
| 4 | 50 | ReLU | 7.522e-07 | 1977 |
| 4 | 50 | Tanh | **4.816e-07** | 2186 |
| 4 | 100 | ReLU | 6.395e-05 | 2361 |
| 4 | 100 | Tanh | 2.861e-05 | 2643 |

The neural network with the minimum MSE has 4 hidden layers, each with 50 neurons. The activation function of its hidden layers is Tanh. Additionally, some conclusions can be made from Table 4. Firstly, the networks with ReLU as the activation function take less time for training. Secondly, the networks with Tanh as the activation function achieve smaller MSEs. Thirdly, the network with 4 hidden layers and 100 neurons in each hidden layer has overfitted during the training process.

The variation of MSE of the neural network with 4 hidden layers and with 50 neurons for each layer is shown in Fig. 6. MSE finally converges to 4.816e-07, which transforms into the mean absolute error (MAE) of the DNN's output: [0.004241 km/s; 0.000232 rad; 0.000152 rad].

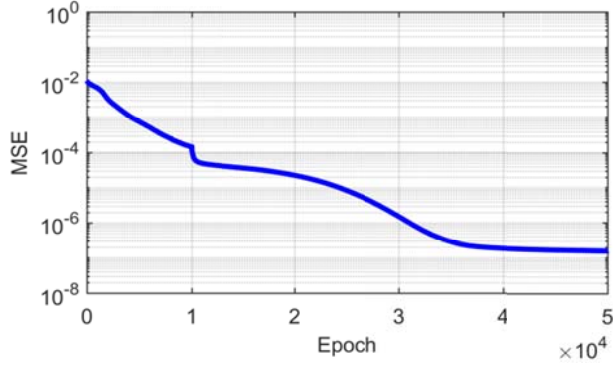

**Fig. 6 MSE of the selected DNN during the training process.**

In order to verify the prediction accuracy of the trained DNN, 1000 new samples that are different from the trained samples were randomly regenerated with the algorithm in part B of Section III to examine the performance of the trained DNN. The initial velocity $v_0$, which is the exact solution of the J2-perturbed Lambert problem and terminal position $r_f$ are used as reference values. The errors of the Keplerian solutions ($v_d$, $r_{fd}$) and the approximation of the trained DNN ($v_c$, $r_{fc}$) are calculated as follows

$$\begin{cases} \Delta v_{0d} = v_0 - v_d, \Delta r_{fd} = r_f - r_{fd} \\ \Delta v_{0c} = v_0 - v_c, \Delta r_{fc} = r_f - r_{fc} \end{cases} \quad (12)$$

Fig. 7 and Fig. 8 show the comparison between the Keplerian solutions and the approximation of the trained DNN. [$\Delta v_{0dx}$; $\Delta v_{0dy}$; $\Delta v_{0dz}$] and [$\Delta r_{fdx}$; $\Delta r_{fdy}$; $\Delta r_{fdz}$] are the errors of the initial velocity and the terminal position of the Keplerian solutions, respectively. [$\Delta v_{0cx}$; $\Delta v_{0cy}$; $\Delta v_{0cz}$] and [$\Delta r_{fcx}$; $\Delta r_{fcy}$; $\Delta r_{fcz}$] are the errors of the initial velocity and the terminal position after the DNN's corrections, respectively. It can be seen that the mean of these errors (red points in Fig. 7 and Fig. 8) is much closer to 0 after the DNN's correction. The standard deviation of these errors has also reduced significantly after the correction, which is indicated by the length of the blue bars in Fig. 7 and Fig. 8. After the correction by the DNN, the initial velocity error is limited to 10 m/s, and the terminal position error does not exceed 100 km. This proves that the application of the DNN has significantly improved the accuracy of the initial value with respect to a simple Keplerian Lambert solution.

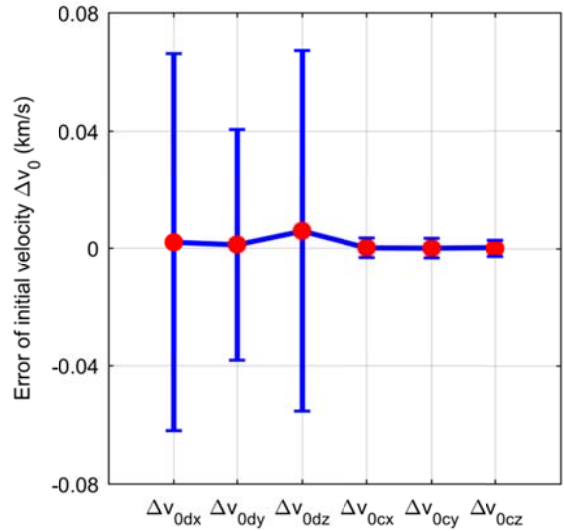

**Fig. 7 The statistical results of the initial velocity errors of the Keplerian solution and the DNN's correction**

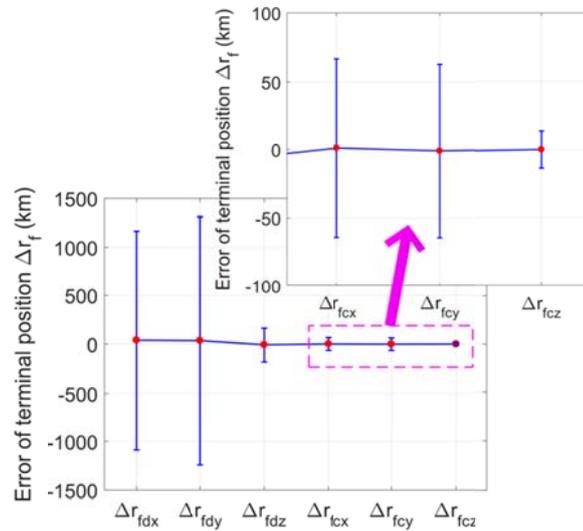

**Fig. 8 The statistical results of the terminal position errors of the Keplerian solution and the DNN's correction**

B. **Performance Analysis for MRPLP**

In this section the proposed DNN-based method is compared against other two methods: a traditional shooting method using Newton's iteration algorithm (SN) and the homotopic perturbed Lambert algorithm (HL) in [15].

When applying the HL, the C++ version of Vinit6 algorithm in literature [27] is employed to implement the HL method in Ref. [15] and to decrease the CPU computation time of HL. The HL is running in Matlab and the MEX function calls the Vinit6 algorithm that is running in visual studio 2015 C++ compiler to analytically propagate the perturbed trajectory. The accuracy tolerance of Vinit6 algorithm is set at $1\times10^{-12}$. The homotopy parameter is defined as the deviation in the terminal position and other details of implementation and settings are the same as these given in Ref. [15]. For the SN and the proposed method, their dynamical models only include the J2 perturbation. For the Vinit6 algorithm, the dynamical model includes the J2, J3 and partial J4 perturbations. However, the magnitudes of J3 and J4 of Jupiter are much smaller than that of J2. Their perturbation effects are very weak compared with that of J2. Therefore, the slight difference in the dynamical model has very limited impact on the number of iterations and running time of the HL since the Vinit6 algorithm has high computational efficiency. Therefore, the comparison among the three methods is still valid.

The performance of the three methods is compared over 11 datasets one per number of full revolutions from 0 to 10. Each dataset has 1000 samples, which are regenerated with the method described in Section III to validate the DNN. The maximum iterations and tolerances of the three methods are listed in Table 5.

Table 5 The maximum iterations and tolerance of three methods

| Algorithm | Tolerance (km) | Maximum iterations |
|---|---|---|
| SN | 0.001 | 2000 |
| HL | 0.001 | 10000 |
| DNN-based method | 0.001 | 2000 |

If the algorithm converges to a solution that meets the specified tolerance within the set number of iterations, it is recorded as a valid convergence, otherwise, as a failed convergence. The result is displayed in Fig. 9 and Fig. 10, in terms of convergence ratio (number of converged solutions over number of samples) and average number of iterations to converge.

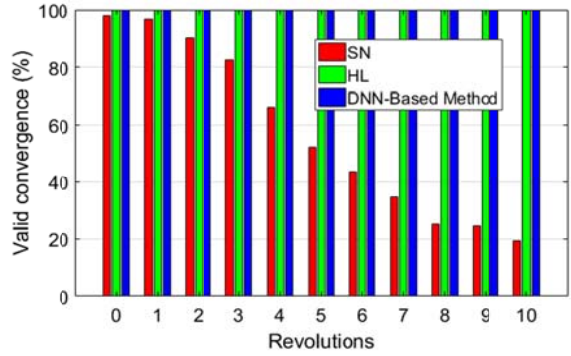

**Fig. 9 The convergence ratio of different algorithms for the Jupiter J2-perturbed Lambert problem**

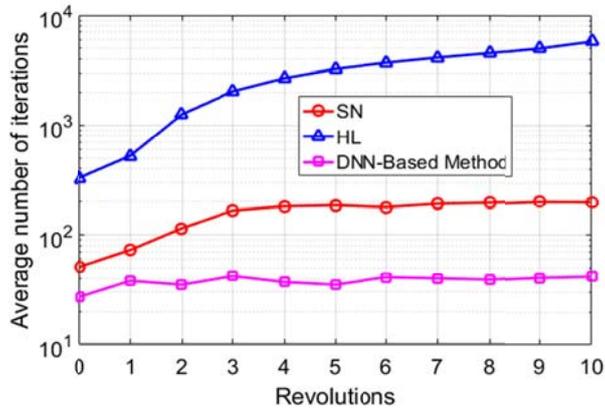

**Fig. 10 Average number of iterations of different algorithms on the Jupiter J2-perturbed Lambert problem**

According to Fig. 9, the HL and the proposed method could converge to the required accuracy in all cases, while the valid convergence ratio of the SN decreases as the number of revolutions increases. Then, according to Fig. 10, the number of iterations of HL appears to increase linearly in log-scale as the number of revolutions increases while the number of iterations of SN and the proposed DNN-based method remain nearly constant. The proposed method requires the least number of iterations. The lack of convergence of the SN with the increase in the number of revolutions is due the growing difference between the exact solution and the solution of the Keplerian Lambert problem. For the same reason the HL progressively requires more iterations to converge. The proposed method mitigates this problem by providing a good initial guess for every number of revolutions.

The average CPU computational time of the three methods is given in Fig. 11, in which the proposed method only accounts the time of the SCM. For zero-revolution case, the average CPU computation time of SN, DNN-Based method and HL are 0.051 seconds, 0.027 seconds and 0.329 seconds, respectively. It is seen that the CPU calculation time of the proposed method is the shortest. This advantage becomes more obvious as the number of the revolution increases because the accurate initial guess obtained using IGG reduces the number of iterations of the SCM. In general, the computational time increases with the increase in the number of revolutions, due to both the increase in the number of iterations and the longer propagation time. As shown in Fig. 11, the computational time of the SN and the proposed method appear to increase linearly with the number of revolutions, while the computational time of HL appears to increase more rapidly. The figure shows that the initial guess obtained with the DNN effectively reduces the number of iterations and provides, as a result, a slower increase of the computational time with the number of revolutions. The computational time of the proposed method is below 0.5 seconds, for the number of revolutions tested in this paper. The average computational time per iteration of SN, HL, and the proposed method are respectively 0.0082 s, 0.0018 s, and 0.0078 s. The proposed method and SN use the same shooting algorithm, for which each iteration needs additional three integral operations to calculate the Jacobian matrix. Their computational time per iteration is higher than that of the HL. However, though the single-iteration of HL takes less time, the HL requires much more iterations than the other two methods, as shown in Fig. 11.

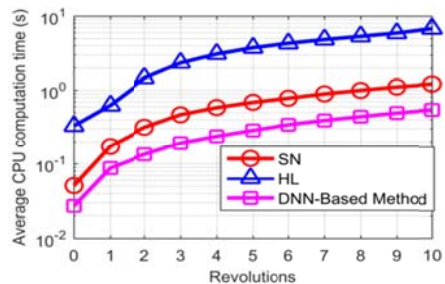

**Fig. 11 Average CPU computational time of different methods for the Jupiter J2-perturbed Lambert problem**

## C. Monte Carlo Analysis

In this section we simulate the repeated use of the DNN-based method by taking a random set of boundary conditions and transfer times and computing multiple J2-perturbed Lambert solutions. Since it is essential to generate samples and train DNN before using the proposed method, the total computational time should include the time of sample generation, the training of the DNN and the SCM. To compare the total CPU time of the above three methods, four sets of Monte Carlo simulations with 1000, 5000, 10000, and 100000 sets of boundary conditions and transfer times are performed. For each set, the number revolutions are equally distributed between 0 and 10. The DNN is trained only once, using 200000 samples and the parameters setting presented in previous section, and is called one time per MC simulation to generate the first guess. The training of DNN was implemented in Python while the solutions of the J2-perturbed Lambert problem using the proposed method, HL and SN run in Matlab. All computations are performed on the personal computer with Intel Core-i7 4.2 GHz CPU and 128GB of RAM. The final results are given in Fig. 12. It can be seen that the efficiency of the proposed method improves with the increase in the number of Lambert solutions to be computed. In particular, when the number of simulations is equal or larger than 5000, the proposed method outperforms the other two methods even when including the cost of the sample generation and the training of the DNN.

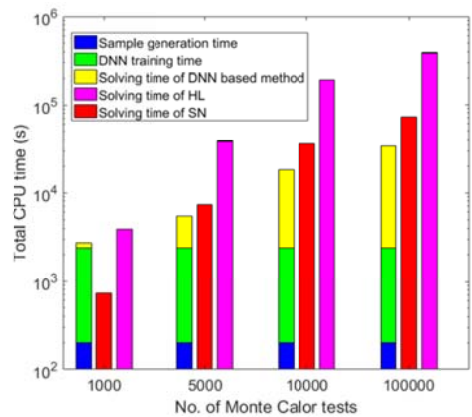

**Fig. 12 Total CPU time of different methods for the Jupiter J2-perturbed Lambert problem**

In addition, two stress cases, where the angle between the initial and terminal positions is 180 deg or 360 deg, have been tested with the proposed method. For each revolution, 100 MC tests are performed for each case. All tests converge successfully and their average CPU computational time is given in Fig. 13, which is similar to the trend in Fig. 11. For the zero revolution case, the CPU computation time of the 180 degree and the 360 degree scenarios are 0.024 seconds and 0.029 seconds, respectively. The case of 360 deg costs a bit more time than the case of 180 deg due to its longer time of flight for each revolution.

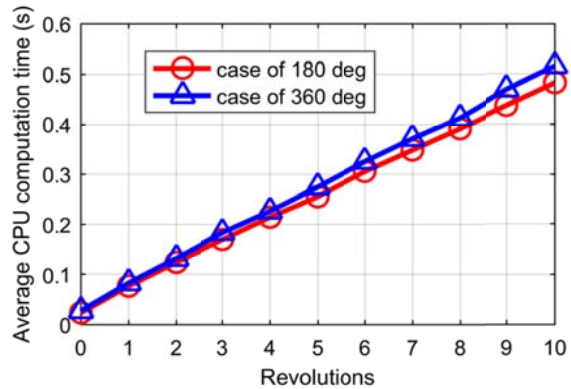

**Fig. 13 Average CPU computational time of two stress cases for the Jupiter J2-perturbed Lambert problem**

## VI.  Conclusion

A fast and novel method using DNN and the finite-difference-based shooting algorithm has been proposed to solve the J2-perturbed Lambert problem. DNN composed of several layers is the extension of conventional artificial neural networks, which has an excellent performance on approximating nonlinear system. The major contribution of the novel method is to use a DNN to generate a first guess of the correction of the initial velocity to solve a J2-perturbed Lambert problem. We demonstrated that the DNN is capable of correcting the initial velocity error of the Keplerian solution and provide good initial values for the subsequent differential correction method. When applied to the Jupiter J2-perturbed Lambert problem, the errors in the initial velocity and terminal position are limited to 5m/s and 100 km, respectively. In addition, when compared to a direct application of a shooting method using Newton's iterations and to a homotopy perturbed Lambert algorithm, the proposed method displayed a computational time that appears to increase linearly with a slope of 0.047 with the number of revolutions. In the application scenario presented in this paper the computational time is less than 0.5 seconds even for ten revolutions. It was also shown that compared to a direct application of a shooting method it provides convergence to the required accuracy in all the cases analyzed in this paper. Thus, we can conclude that the proposed DNN-based generation of a first guess is a promising method to increase robustness and reduce computational cost of shooting methods for the solution of the J2-pertubed Lambert problem.

The method proposed in this paper can be used to solve the J2-perturbed Lambert problem for other celestial bodies, by training the corresponding DNN with the corresponding $J_2$ and $\mu$ parameters. Thus a library of pre-trained DNN could be easily used to have a more general application to missions around any celestial body. On the other hand, adding these dynamical parameters as part of the training set would allow a single more general DNN to be used with all celestial bodies. This latter option is the object of the current investigation.

## Acknowledgments

The work described in this paper was supported by the National Natural Science Foundation of China (Grant No. 11672126), sponsored by Qing Lan Project, Science and Technology on Space Intelligent Control Laboratory (Grant No. 6142208200203 and HTKJ2020KL502019), the Funding for Outstanding Doctoral Dissertation in NUAA (Grant No. BCXJ19-12), State Scholarship from China Scholarship Council (Grant No. 201906830066). The authors fully appreciate their financial supports.

# References

[1] Engels R C, and Junkins J L., "The gravity-perturbed Lambert problem: A KS variation of parameters approach," *Celestial mechanics*, Vol. 24, No.1, 1981, pp. 3-21.
doi: 10.1007/BF01228790

[2] He B, Shen H., "Solution set calculation of the Sun-perturbed optimal two-impulse trans-lunar orbits using continuation theory," *Astrodynamics*, Vol. 4, No. 1, 2020, pp. 75-86.
doi: 10.1007/s42064-020-0069-6

[3] Izzo D, "Revisiting Lambert's problem," *Celestial Mechanics and Dynamical Astronomy*, Vol. 121, No. 1, 2015, pp. 1-15.
doi: 10.1007/s10569-014-9587-y

[4] Bombardelli C, Gonzalo J L, and Roa J., "Approximate analytical solution of the multiple revolution Lambert's targeting problem," *Journal of Guidance, Control, and Dynamics*, Vol. 41, No. 3, 2018, pp. 792-801.
doi: 10.2514/1.G002887

[5] Russell R P., "On the solution to every Lambert problem," *Celestial Mechanics and Dynamical Astronomy*, Vol. 131, No. 11, 2019, pp. 1-33.
doi: 10.1007/s10569-019-9927-z

[6] Der G J., "The superior Lambert algorithm," *Proceedings of the Advanced Maui Optical and Space Surveillance Technologies Conference*, Maui Economic Development Board, Maui, 2011, pp. 462–490.

[7] Armellin R, Gondelach D, and San Juan J F., "Multiple revolution perturbed Lambert problem solvers," *Journal of Guidance, Control, and Dynamics*, Vol. 41, No. 9, 2018, pp. 2019-2032.
doi: 10.2514/1.G003531

[8] Kraige L G, Junkins J L, and Ziems L D., "Regularized Integration of Gravity-Perturbed Trajectories-A Numerical Efficiency Study," *Journal of Spacecraft and Rockets*, Vol. 19, No. 4, 1982, pp. 291-293.
doi: 10.2514/3.62255

[9] Junkins J L and Schaub H., *Analytical mechanics of space systems*, 2nd ed., AIAA, Reston, VA, 2009, pp. 557-561.
doi: 10.2514/4.867231

[10] Arora N, Russell R P, and Strange N, and Ottesen, D., "Partial derivatives of the solution to the Lambert boundary value problem," *Journal of Guidance, Control, and Dynamics*, Vol. 38, No. 9, 2015, pp. 1563-1572.
doi: 10.2514/1.G001030

[11] Woollands R M, Bani Younes A, and Junkins J L., "New solutions for the perturbed lambert problem using regularization and picard iteration," *Journal of Guidance, Control, and Dynamics*, Vol. 38, No. 9, 2015, pp. 1548-1562.
doi: 10.2514/1.G001028

[12] Godal T. "Method for determining the initial velocity vector corresponding to a given time of free flight transfer between given points in a simple gravitational field," *Astronautik*, Vol. 2, 1961, pp. 183-186.

[13] Woollands R M, Read J L, Probe A B, and Junkins J. L., "Multiple revolution solutions for the perturbed lambert problem using the method of particular solutions and picard iteration," *The Journal of the Astronautical Sciences*, Vol. 64, No. 4, 2017, pp. 361-378.
doi: 10.1007/s40295-017-0116-6

[14] Alhulayil M, Younes A B, and Turner J D. "Higher order algorithm for solving lambert's problem," T*he Journal of the Astronautical Sciences*, Vol. 65, No. 4, 2018, pp. 400-422.
doi: 10.1007/s40295-018-0137-9

[15] Yang Z, Luo Y Z, Zhang J, and Tang G J, "Homotopic perturbed Lambert algorithm for long-duration rendezvous optimization," *Journal of Guidance, Control, and Dynamics*, Vol. 38, No. 11, 2015, pp. 2215-2223.
   doi: 10.2514/1.G001198

[16] Li H, Chen S, Izzo D, and Baoying H, "Deep networks as approximators of optimal low-thrust and multi-impulse cost in multitarget missions," *Acta Astronautica*, Vol. 166, 2020, pp. 469-481.
   doi: 10.1016/j.actaastro.2019.09.023

[17] Rubinsztejn A, Sood R, and Laipert F E., "Neural network optimal control in astrodynamics: Application to the missed thrust problem," *Acta Astronautica*, Vol. 176, 2020, pp.192-203.
   doi: 10.1016/j.actaastro.2020.05.027

[18] Izzo D, Märtens M, and Pan B., "A survey on artificial intelligence trends in spacecraft guidance dynamics and control," *Astrodynamics*, 2018, pp. 1-13.
   doi: 10.1007/s42064-018-0053-6

[19] Sánchez-Sánchez C and Izzo D., "Real-time optimal control via Deep Neural Networks: study on landing problems," *Journal of Guidance, Control, and Dynamics*, Vol. 41, No. 5, 2018, pp. 1122-1135.
   doi: 10.2514/1.G002357

[20] Zhu Y and Luo Y Z., "Fast Evaluation of Low-Thrust Transfers via Multilayer Perceptions," *Journal of Guidance, Control, and Dynamics*, Vol. 42, No. 12, 2019, pp. 2627-2637.
   doi: 10.2514/1.G004080

[21] Song Y and Gong S., "Solar-sail trajectory design for multiple near-Earth asteroid exploration based on deep neural networks," *Aerospace Science and Technology*, Vol. 91, 2019, pp. 28-40.
   doi: 10.1016/j.ast.2019.04.056

[22] Cheng L, Wang Z, Jiang F, and Zhou C., "Real-time optimal control for spacecraft orbit transfer via multiscale deep neural networks," *IEEE Transactions on Aerospace and Electronic Systems*, Vol. 55, No. 5, 2018, pp. 2436-2450.
   doi: 10.1109/TAES.2018.2889571

[23] Battin R H. An Introduction to the Mathematics and Methods of Astrodynamics, revised ed., AIAA, VA, 1999, Chap. 6.
   doi: 10.2514/4.861543

[24] Ely T A., "Transforming mean and osculating elements using numerical methods," *The Journal of the Astronautical Sciences*, Vol. 62, No. 1, 2015, pp: 21-43.
   doi: 10.1007/s40295-015-0036-2

[25] Kingma D P, Ba J., "Adam: A method for stochastic optimization," arXiv preprint arXiv:1412.6980, 2014.

[26] Menon A, Mehrotra K, Mohan C K, et al., "Characterization of a class of sigmoid functions with applications to neural networks," *Neural Networks*, Vol. 9, No. 5, 1996, pp: 819-835.
   doi: 10.1016/0893-6080(95)00107-7

[27] Vinti, J. P., Orbital and Celestial Mechanics, Vol. 177, Progress in Astronautics and Aeronautics, AIAA, Reston, VA, 1998, pp. 367–385.
   doi:10.2514/4.866487